\documentclass[12pt]{article}
\usepackage[T2A]{fontenc} 
\usepackage{color}
\usepackage{amsfonts,amssymb}
\newtheorem {theorem}{Theorem}
\newtheorem{lemma}{Lemma}

\newtheorem{rem}{Remark}
\newtheorem{prop}{Proposition}

\begin{document}

\begin{center}
{\Large \bf Codimensions of identities of solvable Lie superalgebras}
\vskip 0.3in
{\bf M.V. Zaicev}\footnote{{\it Zaicev Mikhail Vladimirovich} -- Dr. Sci., Prof., Chair of
Higher Algebra, Faculty of Mech. and Math., MSU,
e-mail: zaicevmv@mail.ru}
{\bf and}
{\bf  D.D. Repov\v s}\footnote{{\it Repov\v s Du\v san Du\v sanovi\v c} -- Dr. Sci., Prof., Chair of Topology and Geometry,
Faculty of Math. and Phys., UL., 
e-mail: dusan.repovs@guest.arnes.si}
\vskip 0.1in

\end{center}
\vskip 0.2in
{\bf Abstract.} We study identities of Lie superalgebras over a field of characteristic zero. We construct
a series of examples of finite-dimensional solvable Lie superalgebras with a non-nilpotent commutator
subalgebra for which PI-exponent of codimension growth exists and is an integer number. 
\vskip 0.1in
\noindent
{\bf Keywords}: identities, codimensions, Lie superalgebras, PI-exponent. 

\section{Introduction}

We study identities of Lie superalgebras over a field $F$ of characteristic zero.  
The existence of a non-trivial identity of an algebra plays a important role in the study of its
properties and structure. For example, if $A$ is an associative finitely generated PI-algebra then its Gelfand -- Kirillov
dimension $GKdim(A)$ is finite and Jacobson radical $J(A)$ is nilpotent. Moreover, if $A$ is simple then
  $\dim A<\infty$. If $A$ and $B$ are two finite-dimensional simple algebras (not necessarily associative) over an algebraically closed field then they are isomorphic if and only if $A$ and $B$ satisfy the same polynomial identities.

Analysis of numerical invariants is one of the fundamental directions in the study of identity relations.
One of the most significant numerical invariants that characterizes the quantity of identities of algebra $A$ is the
sequence $c_n(A), n=1,2,\ldots,$ called the codimension sequence. In the general case, the sequence $\{c_n(A)\}$ has an 
overexponential growth. For example, if $A$ is a free assocative algebra of countable rank then $c_n(A)=n!$. For a free Lie algebra we have $c_n(A)=(n-1)!$. Even if a Lie algebra $L$ satisfies the sufficiently strong identity $[[x_1,x_2,x_3],[y_1,y_2,y_3]]\equiv 0$
then its codimension sequence $\{c_n(L)\}$ grows like $\sqrt{n!}$ (see \cite{V}). Nevertheless, for a wide class     of algebras, the codimension sequence is exponentially bounded. So, for any associative PI-algebra $A$ there is a constant
$a$ such that  $c_n(A)<a^n$ for all  $n\ge 1$ \cite{R} (see also  \cite{L}). If $A$ is an arbitrary 
finite-dimensional algebra, $\dim A=d$, then $c_n(A)\le d^{n+1}$ (see \cite{BD} or  \cite{GZ1}). If  $L$ is an infinite-dimensional simple Lie algebra of Cartan type or Virasoro algebra then  $c_n(A)<a^n$ \cite{M}. Similar restriction 
holds also for any affine Kac -- Moody  algebra \cite{Z1}. If $L$ is a Lie superalgebra with nilpotent commutator
subalgebra, $(L^2)^{t+1}=0$, then the codimension sequence $\{c_n(L)\}$ grows asymptotically not faster than
 $(2t)^n$ \cite{ZM}. For any Novikov algebra,  $A$ the codimension sequence is also exponentially
 bounded, $c_n(A)\le 4^n$ \cite{Dz}.

In the 1980s  S. Amitsur posed a conjecture that the limit of the sequence $\{\sqrt[n]{c_n(A)}\}$
exists and is a non-negative integer for any  associative PI-algebra $A$.
This conjecture  was confirmed in  \cite{GZ2}, \cite{GZ3} and the limit
\begin{equation}\label{eqv1}
exp(A)= \lim_{n\to \infty}\sqrt[n]{c_n(A)},
\end{equation}
was called the PI-exponent of algebra $A$. Later on, the existence and integrality of the limit (\ref{eqv1}) was proved for 
any finite-dimensional Lie algebra \cite{Z2}, Jordan algebra \cite{GSZ}, and some other algebras. It turned out that in the case of
finite-dimensional associative, Lie or Jordan algebra over an algebraically closed field, that the PI-exponent of $A$ is
equal to $\dim A$ if and only if  $A$ is simple.  
 																						
If $A$ is graded by a group $G$, one can also study $G$-graded identities of $A$ and their numerical invariants.
Graded identities form a more precise charateristic than ordinary identities.
For example, if $G=\mathbb Z_2$ then any multilinear identity of degree $n$ is equivalent to the system of $2^n$
graded identities. Therefore it is resonable in the Lie superalgebra case to consider both graded and non-graded identities. 

  It turned out that in the super Lie case, the situation significantly differs from the ordinary Lie or  associative
 case.
In papers \cite{GZ4}, \cite{RZ}, \cite{RZ1} examples of finite-dimensional Lie superalgebras are given for which
graded and  ordinary PI-exponent exist but they are not integer numbers. It was also shown that PI-exponent of simple Lie superalgebra $L$  can be less than $\dim L$.

In the above mentioned examples,   the finite-dimensional Lie superalgebras are not solvable. So, the natural question arises:
Is it true that graded and non-graded exponents exist  for any finite-dimensional solvable Lie superalgebra $L$?
If commutator subalgebra of $L$ is nilpotent then the answer is affirmative (see \cite{ZM}). 
On the other hand, for
solvable Lie superalgebra $L=L_0\oplus L_1$ with nonzero odd component $L_1$, its ideal $L^2$ can be non-nilpotent. 
In   \cite{RZ2} a series of finite-dimensional solvable Lie superalgebras $S(t), t\ge 2$, with non-nilpotent commutator subalgebras was constructed. It was also shown that $exp(S(2))=exp^{gr}(S(2))=4.$ In the present paper we    prove existence and integrality of graded  PI-exponent for any superalgebra $S(t),t\ge 3$. We also compute
the value of this exponent.
 
All necessary information about polynomial identities and their numerical invariants can be found in monographs
 \cite{B}, \cite{Dr}, \cite{GZBook}. 

\section{Preliminaries}

Let $F$ be a field of  characteristic zero and let $F\{X,Y\}$ be the absolutely free algebra over  $F$ with two
infinite sets of generators $X$ and $Y$. Algebra $F\{X,Y\}$ can be naturally endowed with $\mathbb Z_2$-grading  
$F\{X,Y\}=F\{X,Y\}_0\oplus F\{X,Y\}_1$ if we define all generators from $X$ as even and all from $Y$ as odd. If   $L=L_0\oplus L_1$ is some $\mathbb Z_2$-graded algebra over $F$ then a nonassociative polynomial
$f=f(x_1,\ldots,x_m,y_1,\ldots,y_n)\in F\{X,Y\}$ is called a graded identity of algebra $L$ if
$f=f(a_1,\ldots,a_m,b_1,\ldots,b_n)=0$ for any $a_1,\ldots,a_m\in L_0,b_1,\ldots,b_n\in L_1$. The set of all identities  $Id^{gr}(L)$ forms a graded ideal of  $F\{X,Y\}$, invariant under all endomorphisms of 
$F\{X,Y\}$ preserving grading, that is, it is a T-ideal.

Denote by $P_{k,m}$ the subspace of all multilinear polynomials of degree $n=r+m$
on $x_1,\ldots,x_k\in X,y_1,\ldots,y_m\in Y$. It is well-known that the family of all subspaces
 $P_{r,m}\cap Id^{gr}(L),k,m\ge 1$,  uniquely defines $Id^{gr}(L)$ as a T-ideal. Let also
$$
P_{k,n-k}(L)=\frac{P_{k,n-k}}{P_{k,n-k}\cap Id^{gr}(L)}.
$$
Then the value
$$
c_{k,n-k}(L)=\dim P_{k,n-k}(L)
$$
is called the partial $(k,n-k)$-graded codimension, whereas the value     
$$
c_n^{gr}(L)=\sum_{k=0}^n{n\choose k}c_{k,n-k}(L)
$$
is called the $n$-th graded codimension of  $L$.

As in the non-graded case, the sequence of graded codimensions of a finite-dimensional algebra $L$ is exponentially
bounded \cite{BD}. This implies the existence of the limits
$$
\overline{exp^{gr}}(A)= \limsup_{n\to \infty}\sqrt[n]{c_n^{gr}(A)}, \quad
\underline{exp^{gr}}(A)= \liminf_{n\to \infty}\sqrt[n]{c_n^{gr}(A)},
$$
called the upper and the lower graded PI-exponents of $L$, respectively. If the ordinary limit 
$$
exp^{gr}(L)=\lim_{n\to\infty}\sqrt[n]{c_n^{gr}(L)},
$$
exists then it is called the (ordinary) graded PI-exponent of $L$.

Representation theory of symmetric groups is the main tool in the study of numerical characteristics of polynomial
relations. Permutation group $S_n$  acts naturally on multilinear expressions
$$
\sigma\circ f(z_1,\ldots,z_n)=f(z_{\sigma(1)},\ldots,z_{\sigma(n)}).
$$

We recall some elements of representation theory of permutation groups. All details can be found in  \cite{Dzh}.

Denote by $R=FS_m$ the group algebra of group $S_m$. Recall the construction of minimal left ideals of $R$. Let
$\lambda\vdash m$ be a partition  of $m$, that is, an ordered  set of integers $(\lambda_1,\ldots,\lambda_k)$ 
such that $\lambda_1\ge\cdots\ge \lambda_k>0$, $\lambda_1+\cdots+\lambda_k=m$. To this partition corresponds the
so-called Young diagram, that is, the tableau  consisting of $m$ cells, where $\lambda_1$ cells stay in the first row,
$\lambda_2$ cells stay in the second row, etc. Then the Young tableau $T_\lambda$ is the Young diagram  $D_\lambda$ filled 
up  by integers $1,\ldots,m$.

Given Young tableau $T_\lambda$ in $FS_m$, one can construct two subgroups $R_{T_\lambda}$ and  $C_{T_\lambda}$ in $S_m$. The first one is called the row stabilizer and consists of those $\sigma\in S_m$ 
which move integers only within rows     of $T_\lambda$. The second one is called the column stabilizer and consists of permutations which move numbers
 $1,2,\ldots,m$ only within columns of $T_\lambda$. Given Young tableau $T_\lambda$, one can correspond to it  the element 
\begin{equation}\label{eqo0}
e_{T_\lambda}=\left(\sum_{\sigma\in R_{T_\lambda}}\sigma\right)
\left(\sum_{\tau\in C_{T_\lambda}}(-1)^\tau\tau\right),
\end{equation}
of group ring called the Young symmetrizer. It is well-known that Young symmetrizer is  quasi-idempotent, that  is,
$e_{T_\lambda}^2=\gamma e_{T_\lambda}$ where $\gamma\in\mathbb Q$ is a nonzero scalar. Moreover, the left ideal
$R{e_{T_\lambda}}$ is minimal. Its character is denoted by $\chi_\lambda$. Any irreducible left $R$-module
$M$ is isomorphic to some $R{e_{T_\lambda}}$. In this case, its character $\chi(M)$ is equal to $\chi_\lambda$. Recall also that  $R{e_{T_\lambda}}$ 
and
 $R{e_{T_\mu}}$ are isomorphic as $FS_m$-modules if and only if $\lambda=\mu$.

Any finite-dimensional $S_m$-module $M$ can be decomposed into a direct sum of irreducible components
$M=M_1\oplus\cdots\oplus M_t$. In this case, the expression
\begin{equation}\label{eqo1}   
\chi(M)=\sum_{\lambda\vdash m} m_\lambda\chi_\lambda
\end{equation}
means that among $M_1,\ldots,M_t$ there are exactly $m_\lambda$ summands with the character $\chi_\lambda$. The sum of multiplicities $m_\lambda$ in the decomposition (\ref{eqo1}) (that is, the number $t$) is called the length of the
module  $M$.

When we study identities of $\mathbb Z_2$-graded algebras, we need to use the action of direct product of two
symmetric groups on multilinear components of the direct product of two symmetric groups. 
Group  $S_k\times S_{n-k}$ acts  on the space $P_{k,n-k}$.
The intersection $P_{k,n-k}\cap Id^{gr}(L)$ is invariant under this action for any algebra  $L_0\oplus L_1$. Hence $P_{k,n-k}(L)$ is also a $S_k\times S_{n-k}$-module. Any irreducible $S_k\times S_{n-k}$-module is isomorphic to the tensor product $M\otimes N$ of irreducible $S_k$- and $S_{n-k}$-modules, respectively. The character of this
module is denoted by  $\chi_{\lambda,\mu},$ where $\chi_\lambda=\chi(M)$, $\chi_\mu=\chi(N)$. In these notations, decomposition of $P_{k,n-k}(L)$ into irreducible components has the following form: 
\begin{equation}\label{eqo2}
\chi_{k,n-k}(L)=\chi(P_{k,n-k}(L))=\sum_{{\lambda\vdash k\atop\mu\vdash n-k}}
m_{\lambda,\mu}\chi_{\lambda,\mu},
\end{equation}
where $m_{\lambda,\mu}$ is the multiplicity of $\chi_{\lambda,\mu}$ in the decomposition of $\chi_n(L)$.
It follows that
\begin{equation}\label{eqo3}
c_{k,n-k}(L)=\sum_{{\lambda\vdash k\atop\mu\vdash n-k}}
m_{\lambda,\mu} d_\lambda d_\mu,
\end{equation}
where $d_\lambda,d_\mu$ are dimensions of irreducible $S_k$- and  $S_{n-k}$ representations with  the characters
$\chi_\lambda$ and $\chi_\mu$, respectively. 

There is another important series of numerical invariants for estimating the growth of codimensions. 
Value $l_{k,n-k}(L)$ defined as
$$
l_{k,n-k}(L)=\sum_{{\lambda\vdash k\atop\mu\vdash n-k}}
m_{\lambda,\mu},
$$
is called the partial colength of $L$. Here, $m_{\lambda,\mu}$ is an integer on the right hand side of (\ref{eqo2}).  
The total sum
$$
l^{gr}_n(L)=\sum_{k=0}^nl_{k,n-k}(L)=\sum_{k=0}^n
\sum_{{\lambda\vdash k\atop\mu\vdash n-k}}
m_{\lambda,\mu},
$$
is called the graded colength.
 
An important role is played by the estimate  of colength obtained in \cite{Z3}.
 
\begin{lemma}\cite[Theorem 1]{Z3}\label{l1}
Let $L=L_0\oplus L_1$ be a finite-dimensional $\mathbb Z_2$-graded algebra,
$\dim L=d$. Then
$$
l_n^{gr}(L)\le d(n+1)^{d^2+d+1}.
$$
\end{lemma}

\section{Upper estimates for codimension growth}

In this section we obtain an upper bound estimate for codimension growth of Lie superalgebras close to finite-dimensional.
We shall need a technical statement, related to the choice of generators in $S_k\times S_{n-k}$-submodules 
in $P_{k,n-k}$.

\begin{lemma}\label{l2}
Let $M$ be an irreducible $S_k\times S_{n-k}$-submodule in $P_{k,n-k}$ with the character 
$\chi_{\lambda,\mu}$, $\lambda=(\lambda_1,\ldots,\lambda_p)\vdash k$,
$\mu=(\mu_1,\ldots,\mu_q)\vdash(n-k)$. Then there exist $0\ne f=
f(x_1,\ldots,x_k,y_1,\ldots,y_{n-k})\in M$ and decompositions $\{x_1,\ldots,x_k\} 
=
X_1\cup\ldots\cup X_{\lambda_1}$, $\{y_1,\ldots,y_{n-k}\}
=
Y_1\cup\ldots \cup Y_{\mu_1}$
into disjoint subsets such that $f$ is skew symmetric on each of the subsets 
$X_1,\ldots,X_{\lambda_1},Y_1,\ldots,Y_{\mu_1}$. Here, the cardinality $|X_i|$ of each
   $X_i$,  $1\le i\le\lambda_1$ is equal to the height of the $i$-th column of Young diagram
$D_\lambda$, whereas the  cardinality of each $|Y_j|$, $1 \le j\le\mu_1$ is equal to the height of the $j$-th column  of the diagram  $D_\mu$.
\end{lemma}

Proof. 
By hypotheses,  $M$ is isomorphic to $FS_k e_{T_\lambda}\otimes FS_{n-k} e_{T_\mu}$,
where $\lambda\vdash k, \mu\vdash (n-k)$. In particular, $M$ as  $F[S_k\times S_{n-k}]$-module
is generated by elements of the type $(e_{T_\lambda}\otimes e_{T_\mu})h$, where $h=h(x_1,\ldots,x_k,
y_1,\ldots,y_{n-k})$ is a multilinear polynomial. Denote  $h'=e_{T_\lambda}h$. If $e_{T_\lambda}$ has the form (\ref{eqo0}) then we take
$$
h''(x_1,\ldots,x_k,y_1,\ldots,y_{n-k})=
\left(\sum_{\sigma\in C_{T_\lambda}}(-1)^\sigma\sigma\right) h'.
$$

Let $X_1\subseteq \{x_1,\ldots,x_k\}$ consist of all $x_i$ such that indices $i$ are in the first column of
the tableau  $T_\lambda$, $X_2\subseteq \{x_1,\ldots,x_k\}$ consists of all $x_i$ such that indices $i$ are 
in the second column of $T_\lambda$ and so on. Then  $\{x_1,\ldots,x_k\}=X_1\cup\ldots\cup X_{\lambda_1}$ and $h''$ is skew symmetric on each of the sets $X_1,\ldots,X_{\lambda_1}$. Besides, $h''\ne 0$ since $e_{T_\lambda}^2\ne 0$, and
$$
\left(\sum_{\rho\in R_{T_\lambda}}\rho \right)h''=e_{T_\lambda}h'=e_{T_\lambda}^2h.
$$

Next, we set 
$$
f=\left(\sum_{\tau\in C_{T_\mu}}(-1)^\tau\tau \right)h''
$$
and decompose $\{ y_1,\ldots,y_{n-k} \}$ into the union  $Y_1\cup\ldots\cup Y_{\mu_1}$
according to the distribution of indices $y_i$-th  among the columns of $T_\mu$. 
Then  $f\ne 0$ and $Y_1,\ldots,Y_{\mu_1}$ satisfy all required conditions and the proof of the lemma is
completed. 

Recall that any ideal of a Lie superalgebra is by definition  homogeneous in $\mathbb Z_2$-grading. For an upper bound of codimension growth we need the following observation.

\begin{lemma}\label{l3}
Let $L=L_0\oplus L_1$ be a Lie superalgebra and  $I_0\oplus I_1$  its nilpotent ideal of $L$ of finite
codimension, $I^{m+1}=0$. Let also $d_0=\dim(L_0/I_0), d_1=\dim(L_1/I_1)$.
If $\lambda=(\lambda_1,\ldots,\lambda_p)\vdash k$, $\mu=(\mu_1,\ldots,\mu_q)\vdash(n-k)$  
are two partitions such that  $m_{\lambda,\mu}\ne 0$ in the decomposition (\ref{eqo2}) for  $L$,
then $\lambda_{d_0+1}+\cdots+\lambda_p\le m$ and $\mu_{d_1+1}+\cdots+\mu_q\le m$.
\end{lemma}

Proof. Fix a basis  $u_1,u_2,\ldots$ of $L_0$ such that $u_1,\ldots,u_{d_0}$ are linearly independent modulo
 $I_0$, whereas all remaining $u_i$ lie in $I_0$. Similarly, choose a basis $v_1,v_2,\ldots$ in $L_1$ such that
$v_1,\ldots,v_{d_1}$ are linearly independent modulo $I_1$ and $v_j\in I_1, j>d_1$.

Consider an irreducible $S_k\times S_{n-k}$-submodule in $P_{k,n-k}$ with the character
$\chi_{\lambda,\mu}$ and take in $M$ a generator $f=f(x_1,\ldots,x_k,y_1,\ldots,y_{n-k})$
and distributions $X_1,\ldots, X_{\lambda_1}$, $Y_1,\ldots,Y_{\mu_1}$ constructed in Lemma 
\ref{l2}. Suppose that $\lambda_{d_0+1}+\cdots+\lambda_p\ge m+1$. In order to check whether $f$ is an identity of
 $L$ or not it is sufficient to replace variables with elements of fixed bases of corresponding
parity. Let exactly  $t$ first columns of diagram  $D_\lambda$ have the height strictly greater than $d_0$, that is,
$|X_1|,\ldots,|X_t|> d_0, |X_{t+1}|\le d_0$. If we substitute instead of variables from    one of the
sets  $X_i,1\le i\le t$, more than  $d_0$ basis vectors $u_j$ with $j\le d_0$, then we get zero value of $f$ due to 
skew symmety. Otherwise we need to substitute not less  than
$$
N=(|X_1|-d_0)+\cdots+(|X_t|-d_0)
$$
basis elements from  $I$. But since 
$$
N=\lambda_{d_0+1}+\cdots+\lambda_p\ge m+1,
$$
and $I^{m+1}=0$, we again obtain zero value for $f$. Analogously, 
$f\equiv 0$, provided that $\mu_{d_1+1}+\cdots+\mu_q\ge m+1$. Since the inequality
$m_{\lambda,\mu}\ne 0$ implies that  $f$ is not an identity of  $L$, the proof of Lemma \ref{l3} is completed.

Now we estimate dimensions of irreducible components in the decomposition of $P_{k,n-k}(L)$.

\begin{lemma}\label{l4}
Let $\lambda=(\lambda_1,\ldots,\lambda_p)\vdash n$ be a partition of $n$ such that $p\ge d+1$ and $\lambda_{d+1}+\cdots+\lambda_p\le m$. Then, given, $p$ and $m$, the  inequality 
  $d_\lambda\le n^md^n$ holds.
\end{lemma}

Proof. Consider a partition  $\nu=(\lambda_1,\ldots,\lambda_d)$  of the integer  $n'=\lambda_1+\cdots+\lambda_d$. Then  $n-n'\le m$ and by Lemma 6.2.4 from \cite{GZBook}, $d_\lambda\le n^md_\nu$ and by Corollary 4.4.7 from 
\cite{GZBook}, we have $d_\nu\le d^{n'}$.

 \begin{prop}\label{p1} 
Let $L_0\oplus L_1$ be a finite-dimensional Lie superalgebra $\dim L=d$ and let $I=I_0\oplus I_1$ be a nilpotent  
ideal in $L$, $I^{m+1}=0$, $\dim(L_0/I_0)=d_0$, such that  $\dim(L_1/I_1)=d_1$. Then there exists a polynomial $\varphi(n)$ depending only on $m,d,d_0$ and $d_1$ such that
\begin{equation}\label{equ1}
c_{k,n-k}(L)\le\varphi(n)d_0^k d_1^{n-k} \quad \hbox{for all} \quad  0\le k\le n.
\end{equation}
In particular,
\begin{equation}\label{equ1a}
c_n^{gr}(L)\le\varphi(n)(d_0+d_1)^n.
\end{equation}
\end{prop}

Proof. Consider the expression (\ref{eqo3}) for $c_{k,n-k}(L)$. Since all multiplicities  $m_{\lambda,\mu}$ 
are bounded from above by the value $l_n^{gr}(L)$, then by Lemma \ref{l1} we have
\begin{equation}\label{equ2}
c_{k,n-k}(L)=d(n+1)^{d^2+d+1}
\sum_{{\lambda\vdash k\atop\mu\vdash n-k}} d_\lambda d_\mu.
\end{equation}
Skew symmetry considerations  applied in the proof of Lemma \ref{l3} allow us to claim that the height of 
partitions $\lambda$ and $\mu$ (that is the height of corresponding Young diagram)
does not exceed $d$. Clearly, the number of such partitions is less than  $n^d$. 
Hence, applying Lemma \ref{l3}   and Lemma \ref{l4}, we deduce from (\ref{equ2}) the bound (\ref{equ1}) for some polynomial $\varphi(n)$. Now the inequality (\ref{equ1a}) follows from (\ref{equ1}) and the definition of
graded codimension and we have completed the proof of the proposition.

\section{Lie superalgebras of the series $S(t)$}

In this section we define an infinite series of finite-dimensional solvable Lie superalgebras with nonnilpotent
commutator subalgebra. We will use the following agreements. If $A$ is a Lie superalgebra then we denote the product of
elements of $A$ by an ordinary commutator bracket  $[x,y]$. If  $A$ is an associative algebra
then  $[x,y]=xy-yx$. If  $A=A_0\oplus A_1$ is an associative algebra with  $\mathbb Z_2$-grading and  $x$ and $y$ are homogeneous elements from $A$ then
$$
[x,y]=xy-(-1)^{|x||y|}yx,
$$
where $|x|$ is the parity of  $x$, that is, $0$ or $1$. We agree to omit the brackets in the case of
left-normed arrangement, that is, $[x_1,\ldots,x_{k+1}]=[[x_1,\ldots,x_k],x_{k+1}]$ for all $k\ge 2$.

First, let  $R$ be an arbitrary associative algebra with involution $\ast:R\to R$.
Consider the associative algebra $Q=M_2(R)$,
$$
Q= \left\{ \left(
  \begin{array}{cc}
    A & B \\
    C &D \\
  \end{array}
\right) \mid A,B,C,D\in R\right\},
$$
and endow $Q$ by the $\mathbb Z_2$-grading  $Q=Q_0\oplus Q_1$ by setting
$$
Q_0= \left\{ \left(
  \begin{array}{cc}
    A & 0 \\
    0 &D \\
  \end{array}
\right)\right\},\quad
Q_1= \left\{ \left(
  \begin{array}{cc}
    0 & B \\
    C & 0 \\
  \end{array}
\right)\right\}.
$$

It is well-known that $Q$ with the product $[\cdot,\cdot]$ is a Lie superalgebra. Given an associative algebra $R$
with involution, denote by  $R^+$ and  $R^-$ the subspaces of symmetric and skew elements of $R$, respectively:
$$
R^+=\{x\in R~|~~x^\ast=x\},~~ R^-=\{x\in R~|~~x^\ast=-x\}.
$$
Then the subspace    
\begin{equation}\label{eqc1}
L= \left\{ \left(
  \begin{array}{cc}
    x & y \\
    z & -x^\ast \\
  \end{array}
\right) \mid x\in R,y\in R^+,z\in R^-\right\},
\end{equation}
is also Lie superalgebra with the same product as in  $R,$ where
$$
L_0= \left\{ \left(
  \begin{array}{cc}
    x & 0 \\
    0 & -x^\ast \\
  \end{array}
\right)\right\},\quad
L_1= \left\{ \left(
  \begin{array}{cc}
    0 & y \\
    z & 0 \\
  \end{array}
\right)\right\}.
$$

\begin{rem}
In fact, one of the series of simple Lie superalgebras, namely $p(t)$, is constructed  in this way
(see, for example, \cite{S}).
\end{rem}

\begin{rem}
For the Lie superalgebra $L,$ constructed above the following conditions are equivalent:
\begin{itemize}
\item[1)]
$L$ is solvable,
\item[2)]
$L_0$ is solvable Lie algebra,
\item[3)]
$R$ is Lie solvable,
\item[4)]
maximal semisimple subalgebra of  $R$ is commutative.
\end{itemize}
\end{rem}

Thus, the proposed construction gives us a wide class of finite-dimensional solvable Lie superalgebras
with a non-nilpotent (as a rule) commutator subalgebra. As an example we can take 
 algebra of upper triangular $t\times t$ matrices $R=UT_t(F)$, finite-dimensional incidence algebra or any  associative subalgebra in  $UT_t(F)$. We restrict ourselves to the case $R=UT_t(F)$.

Recall the description of involutions  on $UT_t(F)$. One of them $\circ: R\to R$ is the reflection along
secondary diagonal. That is  $e_{ij}^\circ=e_{t+1-j,t+1-i}$ ($e_{ij}$  are the matrix units). It is defined
for all $t\ge 2$. We will call it  {\it orthogonal}.
Another one if defined only for even $t$. Let   $t=2m$. Set    
$$
D=  \left(
  \begin{array}{cc}
    E & 0 \\
    0 & -E \\
  \end{array}
\right) ,
$$
where $E$ is the identity $m\times m$ matrix. Then the map  $s:R\to R$,
$$
X^s=DX^\circ D^{-1},
$$
is also an involution on $R$. It is said to be {\it symplectic}. If we write  $X$ as
$
X=  \left(
  \begin{array}{cc}
    U & V \\
    0 & W \\
  \end{array}
\right) ,
$
then
$
X^s=  \left(
  \begin{array}{cc}
    U^\circ & -V^\circ \\
    0 & W^\circ \\
  \end{array}
\right) .
$

\begin{prop}\cite[Proposition 2.5]{Kosh}\label{p2}
Any involution on  $UT_t(F)$ is equivalent to $\circ$ or $s$.
\end{prop}

{\bf Definition.} Lie superalgebra $(S(t),\ast), t\ge 2,$ is algebra
(\ref{eqc1}),  where $R=UT_t(F)$ and $\ast=\circ$ or $s$ is the orthogonal or symplectic involution on $R$.

Sometimes we will denote both  $(S(t),\circ)$ and  $(S(t),s)$ just by $S(t)$.
We need multiplication formulas in   $L$:
\begin{equation}\label{eqc2}
\left[
\left(
  \begin{array}{cc}
    A & 0 \\
    0 & -A^\ast \\
  \end{array}
\right),
\left(
  \begin{array}{cc}
    0 & B \\
    0 & 0 \\
  \end{array}
\right)
\right]=
\left(
  \begin{array}{cc}
    0 & AB+BA^\ast \\
    0 & 0 \\
  \end{array}
\right),
\end{equation}

\begin{equation}\label{eqc3}
\left[
\left(
  \begin{array}{cc}
    A & 0 \\
    0 & -A^\ast \\
  \end{array}
\right),
\left(
  \begin{array}{cc}
    0 & 0 \\
    C & 0 \\
  \end{array}
\right)
\right]=
\left(
  \begin{array}{cc}
    0 & 0 \\
    -A^\ast C-CA & 0 \\
  \end{array}
\right),
\end{equation}

\begin{equation}\label{eqc4}
\left[
\left(
  \begin{array}{cc}
    A & 0 \\
    0 & -A^\ast \\
  \end{array}
\right),
\left(
  \begin{array}{cc}
    B & 0 \\
    0 & -B^\ast \\
  \end{array}
\right)
\right]=
\left(
  \begin{array}{cc}
    AB-BA & 0 \\
    0 & -(AB-BA)^\ast \\
  \end{array}
\right),
\end{equation}

\begin{equation}\label{eqc5}
\left[
\left(
  \begin{array}{cc}
    0 & B \\
    0 & 0 \\
  \end{array}
\right),
\left(
  \begin{array}{cc}
    0 & 0 \\
    C & 0 \\
  \end{array}
\right)
\right]=
\left(
  \begin{array}{cc}
    BC & 0 \\
    0 & CB \\
  \end{array}
\right).
\end{equation}

Let us introduce few more notations. First note that both involutions $\circ$ and $s$ act similarly on
diagonal matrix units: $e_{ii}^\ast=e_{t+1-i,t+1-i}$. Now, for even $t=2m\ge 2$ or for odd $t=2m+1\ge 3$ we denote
$$
X_i=
\left(
  \begin{array}{cc}
    e_{ii}-e_{ii}^\ast & 0 \\
    0 & e_{ii}-e_{ii}^\ast \\
  \end{array}
\right),
Y_i=
\left(
  \begin{array}{cc}
    0 & e_{ii}+e_{ii}^\ast \\
    0 & 0 \\
  \end{array}
\right),
Z_i=
\left(
  \begin{array}{cc}
    0 & 0 \\
    e_{ii}-e_{ii}^\ast & 0 \\
  \end{array}
\right)
$$
for all  $i=1,\ldots,m,$ and    
$$
E_{ij}=
 \left(
   \begin{array}{cc}
     e_{ij} & 0 \\
      0 & -e_{ij}^\ast \\
   \end{array}
 \right),~1\le i<j\le t,~
I=
\left(
   \begin{array}{cc}
     E & 0 \\
      0 & -E \\
   \end{array}
 \right),~
 Y_0=
\left(
   \begin{array}{cc}
     0 & E \\
      0 & 0 \\
   \end{array}
 \right),
$$
where $E$ is the identity $t\times t$ matrix. The following relations follow from multiplication formulas (\ref{eqc2}) - (\ref{eqc5}):
\begin{equation}\label{eqc5a}
[X_i.Y_j]=[X_i,Z_j]=[X_i,X_j]=0, [Y_i,Z_j]=\delta_{ij}Z_i,~
1\le i,j\le m,
\end{equation}
where $\delta_{ij}$ is the Kronecker symbol and also 
$$
[E_{ik},E_{kj}]=E_{ij},1\le i<k<j\le 2m, [E_{k,k+1},X_{k+1}]=-[E_{k,k+1},X_{k}]=E_{k,k+1}
$$
\begin{equation}\label{eqc6}
[E_{k,k+1},X_{j}]=0, j\ne k,k+1,~[I,E_{ij}]=0, [I,Y_0]=2Y_0.
\end{equation}

At the end of this section we give a lower bound for PI-exponent.

\begin{prop}\label{p3}
Let   $L$ be a Lie superalgebrra of the type $S(t)$. Then $\overline{exp}^{gr}(L)\le 2t$ for even $t$ or $\overline{exp}^{gr}(L)\le 2t-1$ for odd $t$.
\end{prop}

Proof. First note that besides $\mathbb Z_2$-grading, algebra $L$ is also  endowed  by $\mathbb Z$-grading of the
type $L=L^{(0)}\oplus\cdots \oplus L^{(t-1)}$. Initial algebra $R$ has  $\mathbb Z$-grading
$R=R^{(0)}\oplus\cdots \oplus R^{(t-1)}$, where
$$
R^{(k)}=Span\{e_{ij}|~j-i=k\}.
$$
Now, if we put 
$$
L^{(k)}= 
\left\{ \left(
  \begin{array}{cc}
    A & B \\
    C & -A^\ast \\
  \end{array}
\right) \mid A\in R^{(k)},B\in R^+\cap R^{(k)},C\in R^-\cap R^{(k)}\right\},
$$
then multiplication rules (\ref{eqc2}) - (\ref{eqc5}) show that 
$L=L^{(0)}\oplus\cdots \oplus L^{(t-1)}$ is the required $\mathbb Z$-decomposition.
All subspaces $L^{(j)}$ are homogeneous in $\mathbb Z_2$-grading, hence 
$L^{(1)}\oplus\cdots \oplus L^{(t-1)}$ is an ideal of $L$ of codimension $2t$. Since this ideal is nilpotent,  Proposition \ref{p1} completes the proof for even  $t$.

Now let $t=2m+1$. In order to apply Proposition  \ref{p1} again it is enough  to  show that
 $I=<b>+L^{(1)}\oplus\cdots \oplus L^{(t-1)}$ is a nilpotent ideal of  $L$, where
$$
b= 
\left(
  \begin{array}{cc}
    0 & E_{m+1,m+1} \\
    0 & 0 \\
  \end{array}
\right)\in L_1\cap L^{(0)}.
$$
First we need to check that  $[a,b]\in I$, if $a$ is even or odd element from  $L^{(0)}$. If $a$ is even then  $[a,b]=\alpha b, \alpha\in F$, as follows from (\ref{eqc2}) and definition $L^{(0)}$. If
\begin{equation}\label{eqc7}
a= 
\left(
  \begin{array}{cc}
    0 & c \\
    0 & 0 \\
  \end{array}
\right),
\end{equation}
then $[ab]=0$ since the product of any  two elements of the type (\ref{eqc7}) is equal to zero.
On the other hand, if 
$$
a= 
\left(
  \begin{array}{cc}
    0 & 0 \\
    c & 0 \\
  \end{array}
\right),
$$
then $c$ is a diagonal matrix with zero entry at the $(m+1)$-th position. Hence $[a,b]=0$ according to (\ref{eqc5}).

Now let us prove  that $I^{4t}=0$. Let $a=[b_1,\ldots,b_{4t}]$ be a left-normed commutator of homogeneous
both in  $\mathbb Z_2$-  and in $\mathbb Z$-grading elements. If among  $b_i'$s appear at least $t$
 factors from $L^{(1)}\oplus\cdots \oplus L^{(t-1)}$  then  $a=0$, as follows from $\mathbb Z$-grading arguments. On the     other hand, if such factors appear less than $t$ times then $b$ appears at least three times in a row since $I^{(0)}=<b>$. 
In this case, also  $a=0$ since  $(ad~b)^3=0$, where $ad~x$ is the operator of the right hand side
multiplication on  $x$.
Since  
$$
\dim(L/I)=2m+1+2m=2t-1,
$$
 our statement again follows from Proposition \ref{p1} and we have completed the proof.

\section{Exponents of superalgebras of series $S(t)$}

For lower bound of codimension growth we need to consider multialternating polynomials. It will be  
convenient to use the following agreement. If some expression depends on skew symmetric set of arguments
then instead of alternating sum  we will mark these arguments from above by some common symbol - line, tilde, etc.
For example,
$$
\widetilde{x_1}\cdots\widetilde{x_n} =\sum_{\sigma\in S_n}x_{\sigma(1)}\cdots x_{\sigma(n)} 
$$
is the standard polynomial in an associative algebra,
$$
x\overline y z\overline t = xyzt-xtzy,
$$
and
$$
\overline x\overline{\overline y}z \overline{\overline x}\overline y =
x\overline{\overline y}-y\overline{\overline y}z\overline{\overline x}x=
xyzxy-xxzyy-yyzxx+yxzyx.
$$
First consider superalgebras  $S(t)$ with  even $t$. 
\begin{lemma}\label{l5}
Let $S(t)$ be a Lie superalgebra defined by orthogonal or symplectic involution
 $\ast$
  and $t=2m$. Then  $\underline{exp}^{gr}(L)\ge 4m$.
\end{lemma}

Proof. The following relation holds in the algebra of upper triangular matrices.
$$
[[e_{12},\overline e_{11}],\ldots,[e_{n,m+1},\overline e_{mm}]]=
[[e_{12},e_{11}],\ldots,[e_{m,m+1},e_{mm}]]=(-1)^me_{1,m+1}.
$$
It follows that
$$
a_1=
[[E_{12},\overline E_{11}],\ldots,[E_{n,m+1},\overline E_{mm}]]=
$$
\begin{equation}\label{eqe1}
[[E_{12},E_{11}],\ldots,[E_{m,m+1},E_{mm}]]=(-1)^mE_{1,m+1}.
\end{equation}

Expression $a_1$ contains an alternating set of even elements 
$E_{11},,\ldots,E_{mm}$. Let us complicate its construction by adding $m$-alternating odd set. Since
 $[Y_i,Z_i]=X_i, [Y_i,Z_j]=0$ if  $i\ne j$ (see (\ref{eqc5a})) and
  $[E_{k,k+1},X_{k+1}]=[E_{k,k+1},E_{k+1,k+1}]$,~(see (\ref{eqc6})), then
$$
[[E_{12},[\widetilde Y_1,Z_1]],\ldots,[E_{m,m+1},[\widetilde Y_m,Z_{m}]]]=
[[E_{12},[ Y_1,Z_1]],\ldots,[E_{m,m+1},[Y_m,Z_{m}]]]=
$$
$$
=(-1)^m[E_{12},\ldots,E_{m,m+1}]
(-1)^mE_{1,m+1}.
$$

Now we double the number of alternating odd elements. Put 
$$
a_2=[[E_{12},[\widetilde Y_1,Z_1],[Y_1,\widetilde Z_1]],\ldots,[E_{m,m+1},
[\widetilde Y_m,Z_m],[Y_m,\widetilde Z_m]]].
$$
Since $[Z_i,Z_j]=0$ if  $i\ne j$, we can omit an alternation in   $a_2$   not changing the value of
whole expression, that is 
$$
a_2=[[E_{12},[Y_1,Z_1],[Y_1, Z_1]],\ldots,[E_{m,m+1},
[Y_m,Z_m],[Y_m,Z_m]]].
$$
Finally, we put 
$$ 
a_3=[[E_{12},\overline E_{11}, [\widetilde Y_1,Z_1],[Y_1,\widetilde Z_1]],\ldots,[E_{m,m+1},
\overline E_{mm},[\widetilde Y_m,Z_m],[Y_m,\widetilde Z_m]],
$$
$$
[E_{m+1,m+2},\overline E_{m+1,m+1}],\ldots,[E_{2m-1,2m},\overline E_{2m-1,2m-1}],[Y_0,\overline I]].
$$
It follows from multiplication formulas (\ref{eqc5a}), (\ref{eqc6})  that one can omit both alternations in $a_3$ 
preserving the value. In particular,
\begin{eqnarray}\label{eqe2}
a_3=[E_{1,2m-1},[Y_0,I]]=2[E_{1,2m-1},Y_0]=
2
\left(
  \begin{array}{cc}
    0 & e_{1,2m-1}\pm e_{2,2m} \\
    0 & 0 \\
  \end{array}
\right),
\end{eqnarray}
where the sign plus or minus on the right hand side of (\ref{eqe2}) depends on the choice of the involution $\ast$.

Procedure for constructing  element $a_3$ allows us to replicate skew symmetric  sets of even factors $\{E_{i,i},I\}$
as well as odd factors $\{Y_i,Z_i\}$ of $A$. Namely, set
$$ 
A_i^{(0)}=[E_{i,i+1},[Y_i^{(1)},Z_i^{(0)}]],
$$
$$ 
A_i^{(1)}=[A_i^{0},[Y_i^{(2)},Z_i^{(1)}]],
$$ 
$$
-\quad - \quad -
$$
$$ 
A_i^{(p)}=[A_i^{p-1},[Y_i^{(p+1)},Z_i^{(p)}]].
$$ 
Here all  $Y_i^{(j)}$ are copies of element $Y_i$. We use the upper index only for further indication in which
alternation set we will include it. Similar remark holds also for  $Z_i^{(j)}$. 

Further, let
$$ 
A_i^{(p,1)}=[A_i^{(p)},E_{ii}^{(1)}],\ldots,
A_i^{(p,q)}=[A_i^{(p,q-1)},E_{ii}^{(q)}].
$$ 
For $j=m+1,\ldots,2m-1$ we set 
$$
A_j^{(0)}=E_{j,j+1},
$$
$$
A_j^{(1)}=[A_j^{(0)},E_{jj}^{(1)}],
$$
$$
-\quad - \quad -
$$
$$
A_j^{(q)}=[A_j^{(q-1)},E_{jj}^{(q)}].
$$
Finally,
$$
A_{2m}^{(1)}=[Y_0,I^{(1)}],\ldots, A_{2m}^{(q)}=[A_{(2m-1)}^{(q-1)},I^{(q)}].
$$

Let now
$$
W^{(p,q)}=[A_1^{(p,q)},\ldots,A_m^{(p,q)},A_{m+1}^{(q)},\ldots,A_{2m}^{(q)}]
$$
for all  $p,q\ge 1$. Note that for  computing the value of the product $W^{(p,q)}$ it is useful to remember
that the right hand multiplication by $E_{ii}$ commutes with the right hand multiplication by  $[Y_i,Z_i]=X_i$.

Commutator $W{(p,q)}$ depends on
\begin{itemize}
\item[-] $p$ sets of odd elements  $Y_1^{(i)},\ldots,Y_m^{(i)},
         Z_1^{(i)},\ldots,Y_m^{(i)}$, $1\le i \le p$, of the size $2m$;  
\item[-] $q$ sets of even elements  $E_{11}^{(j)},\ldots,E_{2m-1,2m-1}^{(j)},
         I^{(j)}$,  $1\le j \le q$, of the size $2m$;  
\item[-] and also on $4m$ factors  $E_{12},\ldots,E_{2m-1,2m},Y_0,Z_1^{(0)}.
          \ldots,Z_m^{(0)},Y_1^{(p+1)},\ldots,Y_m^{(p+1)}$ outside these sets.
\end{itemize}
Applying to $W^{(p,q)}$ the alternation on the sets of order $2m$, we   get an expression
$$
\widetilde W^{(p,q)}=Alt_1^{(0)}\cdots Alt_q^{(0)}Alt_1^{(1)}\cdots Alt_p^{(1)}(W^{(p,q)}).
$$ 
Here, $Alt_j^{(0)}$ is the alternation on  $E_{11}^{(j)}, E_{2m-1,2m-1}^{(j)}$ 
and
 $I^{(j)}$, whereas
$Alt_i^{(1)}$ is the alternation on $Y_1^{(i)},\ldots,Y_m^{(i)},Z_1^{(i)},\ldots,Z_m^{(i)}$. 

As in computing expressions  $a_1,a_2$ 
and $a_3$, alternation in  $\widetilde W^{(p,q)}$ does not play any role,
that is, 
\begin{equation}\label{eqe3}
\widetilde W^{(p,q)}=W^{(p,q)}=\pm 2^q [E_{1,2m-1},Y_0]\ne 0.
\end{equation}

Now we construct $\widetilde w^{(p,q)}$  in $F\{X,Y\}$ using the same procedure as for 
$\widetilde W^{(p,q)}$, only changing  $E_{12},\ldots,E_{2m-1,2m}$ by even generators  $x_{12},\ldots,x_{2m-1,2m}$,
 $E_{11}^{(j)},\ldots,E_{2m-1,2m-1}^{(j)},I^{(j)}$
by even generators   $x_{1}^{(j)},\ldots,x_{2m,2m}^{(j)}$, $Y_1^{(i)},\ldots,Y_m^{(i)}$
by odd  $y_1^{(i)},\ldots,y_m^{(i)}$,  $Z_1^{(i)},\ldots,Z_m^{(i)}$
by odd $z_1^{(i)},\ldots,z_m^{(i)}$, and $Y_0$ by odd $y_0$.

Element $\widetilde w^{(p,q)}$  includes $q$ skew symmetric  sets of even variables
 $X^{(j)}=\{x_1^{(j)},\ldots,x_{2m}^{(j)}\}, 1\le j\le q$,  and $p$
skew symmetric  sets of odd variables  $Y^{(i)}=
\{y_1^{(i)},\ldots,y_m^{(i)},z_1^{(i)},\ldots,z_m^{(i)}\}$. In addition to  these variables,
$\widetilde w^{(p,q)}$ contains $4m$ variables $x_{12},\ldots$, $x_{2m-1,2m},y_0,
y_1^{(p+1)},\ldots,y_m^{(p+1)}$, $z_1^{(0)},\ldots,z_m^{(0)}$ not participating in alternations.

Fix now  $n=2mp+2mq+4m$ and $k=2mq+2m-1$. Then $n-k=2mp+2m=1$. Subgroup $H=S_{2mq}\times S_{2mp}$ of 
$S_k\times S_{n-k}$ acts on the space $P_{k,n-k}$. Left factor $S_{2mq}$ acts on  
$\overline X=X^{(1)}\cup\cdots\cup X^{(q)}$, whereas $S_{2mp}$ acts on $\overline Y=Y^{(1)}\cup\cdots\cup Y^{(p)}$. Relation (\ref{eqc3}) means that  $\widetilde w (p,q)$ is not an identity of $L$.
Moreover, $\varphi(\widetilde w(p,q))\ne 0$ for the evaluation $\varphi$ such that
$\varphi(\overline  X)\subseteq V_0,\varphi(\overline Y)\subseteq(V_1)$, where
 $V_0=L_0\cap L^{(0)},V_1=L_1\cap L^{(0)}$ are subspaces of dimension $2m$. 
It follows from  the structure of essential idempotent (see (\ref{eqo0})) and skew symmetry of
 $\widetilde w(p,q)$ that in the decomposition of $FH$-submodule in 
$P_{k,n-k}$ generated by  $\widetilde w(p,q)$ appear only irreducible components with the character
 $\chi_{\lambda,\mu}$, where
$$
\lambda=(q^{2m})=(\underbrace{q,\ldots,q}_{2m}),~~
\mu=(p^{2m})=(\underbrace{p,\ldots,p}_{2m})
$$
are two rectangular partitions. From this it follows that $c_{k,n-k}(L)\ge \deg\chi_{\lambda,\mu}
=d_\lambda d_\mu$.

It is well-known that dimension of irreducible representation with rectangular Young diagram
is exponential, where the ratio of exponent is the height of the diagram. For example, by Lemma 5.10.1 
from \cite{GZBook}, for $\nu=s^d\vdash N=sd$ for all $s$ large enough the following inequality holds
$$
d_\nu>N^{-\frac{d(d-1)}{2}} d^N
$$
provided that $d$ is fixed.  In our case for  $k=N+2m-1, N=2mq$, we have 
$$
d_\lambda>\frac{1}{N^{m(m-1)}}(2m)^{k-2m+1}>\frac{1}{n^{m(m-1)}}\cdot\frac{(2m)^k}{(2m)^{2m-1}}.
$$
Similarly, 
$$
d_\mu>\frac{1}{n^{m(m-1)}}\cdot\frac{(2m)^{n-k}}{(2m)^{2m+1}}.
$$
Hence we have proved the inequality 
\begin{equation}\label{eqe4}
c_{k,n-k}>\frac{1}{n^{2m(m-1)}(2m)^{2m}} (2m)^n
\end{equation}
for $k=2mq+2m-1, n-k=2mp+2m+1$.

To obtain analogous lower bound estimate for  $c_{k,n-k}$ for  arbitrary $k$ and $n-k$ large enough, note that
$$
[\widetilde W^{(p,q)},\underbrace{E_{11},\ldots,E_{11}}_r]\ne 0,~~
[\widetilde W^{(p,q)},\underbrace{[Y_1,Z_1],\ldots,[Y_1,Z_1]}_r]\ne 0
$$
in $L$ for any $r\ge 1$. Hence the polynomial 
\begin{equation}\label{eqe5}
[\widetilde w^{(p,q)},x_1,\ldots,x_i,y_1,\ldots,y_j]
\end{equation}
is not a graded identity of $L$.

Now for an arbitrary pair $k,n$ we can find $0\le i,j\le 2m-1, p$  and $q$ such that
$k=k_0+i,n=n_0+j$ where $k_0=2mq+2m-1,n_0-k_0=2mp+2m+1$. Using  the same arguments for
polynomial (\ref{eqe5}) that we apply for $\widetilde w^{(p,q)}$, we obtain lower bound
$$
\dim P_{k,n-k}(L)\ge d_\lambda d_\mu \ge\frac{(2m)^{k_0}}{n_0^{m(m-1)} (2m)^{2m-1} }
\cdot
\frac{(2m)^{n_0-k_0}}{n_0^{m(m-1)} ( 2m )^{2m+1} }\ge
$$
$$
\ge \frac{(2m)^{n_0}  }  {n^{2m(m-1)}(2m)^{4m} }=
\frac{(2m)^{n}  }  {n^{2m(m-1)}(2m)^{4m+i+j} } .
$$
Hence, taking into account (\ref{eqe4}), we get the restriction 
$$
c_{k,n-k}(L)\ge \frac{(2m)^{n}}{n^{2m(m-1) (2m)^{8m}}}
$$
for all  $k$ and $n-k$ large enough. From this follows the inequality
\begin{equation}\label{eqe6}
c_{n}^{gr}(L) \ge \frac{(2m)^{n}}{n^{2m(m-1) (2m)^{8m}}}
\sum_{k=C+1}^{n-C}{n\choose k},
\end{equation}
where $C$ is some constant
 depending only on $m$. Since the sum of binomial 
coefficients is equal to $2^n$, we obtain from (\ref{eqe6}) the estimate
for lower limit,
$$
\underline{exp}^{gr}(L)\ge 4m 
$$
and complete the proof  of our lemma.

Now consider the case of odd  $t$.

\begin{lemma}\label{l6}
Let $t=2m+1$ and  $l=s(t)=(S(t),\circ)$. Then  $\underline{exp}^{gr}(L)\ge 2t-1=4m+1$.
\end{lemma}

Proof. The proof in this case largely repeats the proof of the previous lemma and we will omit   details.
Values $A_i^{(p)}$ and $A_i^{(p,q)}$, $1\le i\le m$ remain the same if we preserve previous notations.
Also  $A_j^{(1)},\ldots,A_j^{(q)}$, $m+1\le j\le 2m-1$ do not change. 
Elements $A_{2m}^{(1)},\ldots, A_{2m}^{(q)}$ 
are defined by induction:
$A_{2m}^{(1)}=E_{2m,2m+1}^{(1)},\ldots,A_{2m}^{(q)}=[A_{2m}^{(q-1)},E_{2m,2m}^{(q)}]$, whereas
$A_{2m+1}^{(q)}$ is defined as $A_{2m}^{(q)}$ in Lemma \ref{l5}. In the expression for $W^{(p,q)}$
we need to replace the last factor $A_{2m}^{(q)}$ by $A_{2m+1}^{(q)}$.
 
Modified element  $\widetilde w^{(p,q)}$ depends on $q$ skew symmetric  sets of even variables of order $2m+1$, 
on $p$ skew symmetric of odd sets variabes of order $2m$ and has total degree $n=2mp+(2m+1)q+4m+1$. Partially modified are lower bounds for $d_\lambda$ and $d_\mu$:
$$
d_\lambda>\frac{ 1 }{n^{m(2m+1)}}\cdot \frac{(2m+1)^k  }{(2m+1)^{2m}  },~~
d_\mu>\frac{ 1 }{n^{m(2m+1)}}\cdot \frac{(2m)^{n-k}  }{(2m)^{2m+1}  },
$$
where $\lambda=(q^{2m+1}),\mu=(p^{2m})$, and also for $c_{k,n-k}(L)$:
$$
c_{k,n-k}(L) \ge \frac{(2m+1)^k(2m)^{n-k}}{n^{2m(2m+1)} (2m+1)^{8(m+1)} }.
$$
Therefore lower bound estimate for graded codimension takes the form:
$$
c_n^{gr}(L) \ge
\frac{ 1 }{n ^{ 2m(2m+1) }( 2m+1 )^{8(m+1)  }} 
\sum_{k=C+1}^{n-C} {n\choose k}(2m+1)^k (2m)^{n-k},
$$
from which we get the following inequality 
$$
\underline{exp}^{gr}(L) \ge 4m+1,
$$
and we have completed the proof.

As an immediate consequence of Lemma  \ref{l5}, Lemma \ref{l6} and Proposition  \ref{p3}, we get the main result
of the paper.

\begin{theorem}\label{t1}
Let  $L=(S(t),\ast)$ be a Lie superalgebra of the type $S(t)$, where $\ast$ is orthogonal or symplectic involution.
Then  the graded PI-exponent of $L$ exists and
\begin{itemize}
\item[$\bullet$]
$\exp^{gr}(L)=2t$ if $t$ is even;
\item[$\bullet$]
$\exp^{gr}(L)=2t-1$ if $t$ is odd.
\end{itemize}
\end{theorem}

\vskip .2in
{\bf Acknowledgements.} The first author was partially supported by the Russian Science Fund, grant  No. 22-11-00052.
The second author was supported by the Slovenian Research and Innovation Agency, grants No. P1-0292, N1-0278, N1-0114, N1-0083,
J1-4031, J1-4001.

\vfill\eject


\begin{thebibliography}{99}
 \bibitem{V}
 I.B. Volichenko.  Varieties of Lie algebras with identity  $[[X_1,X_2,X_3],[X_4,X_5,X_6]]=0$ over a field
 of characteristic zero, {\em Siberian Math. J.,} {\bf 25}:3 (1984), 370 -- 382.
  
\bibitem{R}
A. Regev, Existence of identities in $A\otimes B$, {\em  Israel J. Math.}, \textbf{11}
(1972), 131 -- 152.

\bibitem{L} V.N. Latyshev, On Regev's theorem on identities in a tensor product  of PI-algebras (in Russian), 
{\em Usp. Mat. Nauk},  \textbf{27}:4 (1972),  213 -- 214.

\bibitem{BD}
Yu. Bahturin, V. Drensky, Graded polynomial identities of matrices,
{\em Linear Algebra Appl.}, \textbf{357} (2002), 15 -- 34.

\bibitem{GZ1}
A. Giambruno, M. Zaicev, Codimension growth of special simple Jordan algebras, {\em Trans. Amer. Math. Soc.},
 {\bf 362}:6 (2010). 3107 -- 3123.

\bibitem{M}
S.P. Mishchenko, Growth in varieties of Lie algebras, {\em  Russian Math. Surveys}, {\bf 45} (1990), 27 -- 52.

\bibitem{Z1}
M.V. Zaicev, Varieties of affine Kac-Moody algebras,
 {\em Math. Notes}, {\bf 62}:1 (1997), 80 -- 86.

\bibitem{ZM}
M.V. Zaitsev, S.P. Mishchenko, 
Identities for Lie superalgebras with a nilpotent commutator subalgebra, {\em Algebra and Logic}, {\bf 47}:5 (2008), 348 -- 364.

\bibitem{Dz}
A.S. Dzhumadil'daev,  Codimension growth and non-Koszulity of Novikov operad,
{\em Comm. Algebra}, {\bf 39}:8 (2011), 2943 -- 2952.

\bibitem{GZ2}
A. Giambruno, M. Zaicev, On codimension growth of finitely generated associative algebras, 
{\em Adv. Math.}, {\bf 140}:2 (1998). 145 -- 155. 

\bibitem{GZ3}
A. Giambruno, M. Zaicev, Exponential codimension growth of PI algebras: an exact estimate, 
{\em Adv. Math.}, {\bf 142}:2 (1999), 221 -- 243. 

\bibitem{Z2}
M.V.Zaicev, Integrality of exponents of codimension growth of finite-dimensional Lie algebras,
{\em Izv. Math.}, {\bf 66}:3 (2002), 467 -- 487.

\bibitem{GSZ}
A. Giambruno, I. Shestakov, M. Zaicev, Finite-dimensional non-associative algebras and 
codimension growth,  
{\em Adv. Appl. Math.}, {\bf 47}:1 (2011), 125 -- 139. 

\bibitem{GZ4} 
A. Giambruno, M. Zaicev, On codimension growth of finite-dimensional Lie superalgebras,
{\em J. Lond. Math. Soc. (2)}, {\bf 85}:2 (2012), 534 -- 548.

\bibitem{RZ}
D. Repov\v s, M. Zaicev, Graded identities of some simple Lie superalgebras, {\em Algebr. 
Represent. Theory}, {\bf 17}:5 (2014), 1401 -- 1412. 

\bibitem{RZ1}
D. Repov\v s, M. Zaicev,
Graded codimensions of Lie superalgebra b(2), {\em J. Algebra}, {\bf 422} (2015), 1 -- 10. 

\bibitem{RZ2} 
D. Repov\v s, M. Zaicev, Codimension growth of solvable Lie superalgebras, {\em
J. Lie Theory}, {\bf 28}:4 (2018), 1189 -- 1199.

\bibitem{B}
Yu.A. Bahturin, {\em Identical relations in Lie algebras,} VNU Science Press, b.v., Utrecht, 1987. 

\bibitem{Dr}
V.  Drensky, {\em Free algebras and PI-algebras. Graduate course in algebra},
 Springer-Verlag Singapore, Singapore, 2000. 

\bibitem{GZBook}
A. Giambruno, M. Zaicev, {\em 
 Polynomial Identities and Asymptotic 
Methods}, Mathematical Surveys and Monographs. {\bf 122}. Amer. Math. 
Soc., Providence, RI, 2005.

\bibitem{Dzh}
G.D. James,
The Representation Theory of the Symmetric Groups, 1978 Lecture Notes in Math.682 (Berlin: Springer).

\bibitem{Z3}
M.V. Zaicev, 
Graded identities in finite-dimensional algebras, {\em Moscow University Mathematical
Bulletin}, {\bf 70}:5 (2015), 234 -- 235.

\bibitem{S}
M. Scheunert,
{\em The theory of Lie superalgebras.
An introduction}, Lecture Notes in Mathematics, 716. Springer, Berlin, 1979.

\bibitem{Kosh}
O.M. Di Vincenzo, P. Koshlukov, R. La Scala, Involutions for upper triangular 
matrix algebras, {\em Adv.  Appl. Math.}, {\bf 37}:4 (2006), 541 -- 568.

\end{thebibliography}
\end{document}